\documentclass[12pt]{amsart}
\usepackage{amsmath,amssymb,amsthm,latexsym, xcolor,xspace}
\usepackage{MnSymbol}
\usepackage{amscd}
\usepackage{mathrsfs}
\usepackage{graphicx}
\usepackage{xspace}
\usepackage{color}
\usepackage{enumitem}
\usepackage{hyperref}

\newtheorem{thm}{Theorem}[section]

\newtheorem{cor}[thm]{Corollary}

\theoremstyle{definition}

\newtheorem{example}[thm]{Example}

\theoremstyle{remark}

\newtheorem*{Claim-non}{Claim}

\DeclareMathOperator{\Aut}{{\rm Aut}}

\DeclareMathOperator{\End}{{\rm End}}
\DeclareMathOperator{\QZ}{{\rm QZ}}
\DeclareMathOperator{\PSL}{{\rm PSL}}
\DeclareMathOperator{\LN}{\mathcal{LN}}
\DeclareMathOperator{\LC}{\mathcal{LC}}
\DeclareMathOperator{\rist}{\mathrm{Rist}}
\DeclareMathOperator{\Comm}{\mathrm{Comm}}

\DeclareMathOperator{\con}{\sf{con}}
\DeclareMathOperator{\parb}{\sf{par}}
\DeclareMathOperator{\nub}{\sf{nub}}
\DeclareMathOperator{\lev}{\sf{lev}}

\def\ident{1}

\title[Totally disconnected locally compact groups]{A totally disconnected invitation to locally compact groups}
\author{Pierre-Emmanuel Caprace}
\address{UCLouvain, Belgium}
\author{George A. Willis}
\address{University of Newcastle, Australia}
\date{October 12, 2021} 
\begin{document}

\maketitle

\begin{abstract}
We present a selection of results contributing to a structure theory of totally disconnected locally compact groups. 
\end{abstract}

\section{Introduction}	 
Locally compact groups have attracted sustained attention because, on one hand, rich classes of these groups have fruitful connections with other fields and, on the other, they have a well-developed theory that underpins those connections and delineates group structure. Salient features of this theory are the existence of a left-invariant, or Haar, measure; and the decomposition of a general group into pieces, many of which may be described concretely and in detail. 

Haar measure permits representations of a general locally compact group by operators on spaces of measurable functions, and is thus the foundation for abstract harmonic analysis. Connections with partial differential equations, physics and number theory come about through these representations. Locally compact groups are the largest class for which an invariant measure exists and for which harmonic analysis can be done in this form, as was shown by A.~Weil \cite{Weil}.

The decomposition theory of an arbitrary locally compact group $G$ begins with the short exact sequence 
$$
0 \to G^\circ \to G \to G/G^\circ \to 0,
$$ 
in which the closed normal subgroup $G^\circ$ is the connected component of the identity. The Gleason--Yamabe theorem \cite[Th.~6.0.11]{Tao} applies to $G^\circ$ to show that it is a projective limit of connected Lie groups, and powerful tools from the theory of Lie groups may thus be brought to bear on $G^\circ$. Groups occurring in physics and differential equations are often Lie groups. The quotient $G/G^\circ$ is a totally disconnected locally compact group (abbreviated \textbf{tdlc group}). 
    Lie groups over local fields are important examples of tdlc groups having links to number theory and algebraic geometry (see for example \cite{LafforgueICM,ScholzeICM}).  Unlike the connected case however, many other significant tdlc groups, such as the automorphism groups of locally finite trees first studied in \cite{Tits_arbre}, cannot be approximated by Lie groups. While substantial progress has been made with our understanding of tdlc groups much remains to be done before it could be said that the structure theory has reached maturity. This article surveys our current state of knowledge, much of which is founded on a theorem of van Dantzig, \cite{vD}, which ensures that a tdlc group $G$ has a basis of identity neighbourhoods consisting of compact open subgroups. 

Decompositions of general tdlc groups are described in Section~\ref{sec:dec}. This section includes a discussion of the so-called elementary groups, which are those built from discrete and compact groups by standard operations. Discrete and compact groups are large domains of study in their own right and it is seen how elementary groups can be factored out in the analysis of a general tdlc group. Simple groups are an important aspect of any decomposition theory and what is known about them is summarised in Section~\ref{sec:simple}. This includes a local structure theory and  the extent to which local structure determines the global structure of the group. Section~\ref{sec:scale} treats scale methods, which associate invariants and special subgroups to abelian groups of automorphisms and which in some circumstances substitute for the Lie methods available for connected groups. 
A unifying theme of our approach is the dynamics of the conjugation action: Section~\ref{sec:dec} is concerned with the conjugation action of $G$ on its closed subnormal subgroups, Section~\ref{sec:simple} uses in an essential way   the conjugation action of $G$ on its closed subgroups, especially those that are locally normal, while Section~\ref{sec:scale} concerns the dynamics of the conjugation action of cyclic subgroups (and, more generally, flat subgroups) on the topological space $G$. Section~\ref{sec:future} highlights a few open questions and directions for further research.


\section{Decomposition theory}\label{sec:dec}

\subsection{Normal subgroup structure}

Finite groups, Lie groups and algebraic groups  constitute three of the most important classes of groups. Their respective structure is deep and far-reaching. One of their common themes consists in reducing  problems concerning a given group $G$ in one of these classes to problems about simple groups in the corresponding class, and then tackling the reduced problem by invoking classification results.  Striking illustrations of this approach in the case of finite groups can be consulted in R.~Guralnick's ICM address \cite{Guralnick_ICM}. 

Since the category of locally compact groups contains all discrete groups, hence all groups, developing a similar theory for locally compact groups is hopeless. Nevertheless, the possibility to construct meaningful `decompositions of locally compact groups into simple pieces' has been highlighted in \cite{CaMo}. Wide-ranging results have subsequently been established by C.~Reid and P.~Wesolek in a series of papers \cite{ReidWes_Ess, ReidWes_PLMS}, some of whose results are  summarized below. A more in-depth survey can be consulted in \cite{Reid_survey}.

Given   closed normal subgroups $K, L$  of a locally compact group $G$, the quotient group $K/L$ is called a \textbf{chief factor} of $G$ if $L$ is strictly contained in $K$ and for every closed normal subgroup $N$ of $G$ with $L \leq N \leq K$, we have $N =L$ or $N=K$. Given a closed normal subgroup $N$ of $G$, the quotient $Q = G/N$ is a chief factor if and only if $Q$ is \textbf{topologically simple}, i.e. $Q$ is non-trivial and the only closed normal subgroups of $Q$ are $\{1\}$ and $Q$. More generally, every chief factor $Q=K/L$ is \textbf{topologically characteristically simple}, i.e. the only closed subgroups of $Q$ that are invariant under all homeomorphic automorphisms of $Q$ are $\{1\}$ and $Q$. A topological group is called \textbf{compactly generated} if it has a compact generating set.

\begin{thm}[{See \cite[Th.~1.3]{ReidWes_Ess}}]\label{th:essentially_chief}
Every compactly generated tdlc group $G$ has a finite series  $\{1\} =G_0 < G_1 < G_2 < \dots < G_n = G$ of closed normal subgroups such that  for all $i =1, \dots, n$, the quotient $G_{i}/G_{i-1}$ is compact, or discrete infinite, or a chief factor of $G$ which is  non-compact, non-discrete and second countable.
\end{thm}

A normal series as in Theorem~\ref{th:essentially_chief} is called an \textbf{essentially chief series}. 
That statement is obviously meaningless for compact or discrete groups. Let us illustrate Theorem~\ref{th:essentially_chief} with two examples. 

\begin{example}\label{ex:1}
Let $I$ be a set and for each $i \in I$, let $G_i$ be a tdlc group and $U_i \leq G_i$ be a compact open subgroup. The \textbf{restricted product} of $(G_i, U_i)_{i \in I}$, denoted by $\bigoplus_{i \in I} (G_i, U_i)$, is the subgroup of $\prod_{i \in I} G_i$ consisting of those tuples $(g_i)_{i\in I}$ such that $g_i \in U_i$ for all but finitely many $i \in I$. It is endowed with the unique tdlc group topology such that the inclusion $\prod_{i \in I} U_i \to \bigoplus_{i \in I} (G_i, U_i)$ is continuous and open. Given a prime $p$, set $M(p) = \bigoplus_{n \in \mathbf Z} (\PSL_2(\mathbf Q_p), \PSL_2(\mathbf Z_p))$. The cyclic group $\mathbf Z$ naturally acts on $M(p)$ by shifting the coordinates. The semi-direct product $G(p) = M(p) \rtimes \mathbf Z$ is a compactly generated tdlc group, with an essentially chief series given by $\{1\} < M(p) < G(p)$. The group $M(p)$ is not compactly generated. It has minimal closed normal subgroups, but does not admit any finite essentially chief series, which illustrates the necessity of the compact generation hypothesis in Theorem~\ref{th:essentially_chief}.
\end{example}
\begin{example}\label{ex:2}
A more elaborate construction in \cite[\S 9]{ReidWes_PLMS} yields an example of a compactly generated tdlc group $G'(p)$ with an essentially chief series given by $\{1\} < H(p) < G'(p)$ such that $G'(p)/H(p) \cong \mathbb Z$ and $H(p)$ has  a nested chain of closed normal subgroups $(H(p)_n)$ indexed by $\mathbf Z$, permuted transitively by the conjugation $G'(p)$-action, and  such that $H(p)_n/H(p)_{n-1} \cong M(p)$ for all $n \in \mathbf Z$. 
\end{example}

A tdlc group is compactly generated if and only if it is capable of acting continuously,  properly, with finitely many vertex orbits, by automorphisms on a connected locally finite graph. For a given compactly generated tdlc group $G$, vertex-transitive actions on graphs are afforded by the following construction.   Given a compact open subgroup $U < G$, guaranteed to exist by van Dantzig's theorem, and a symmetric compact generating set $\Sigma$ of $G$, we construct a graph $\Gamma$ whose vertex set is the coset space $G/U$ by declaring that the vertices $gU$ and $hU$ are adjacent if $h^{-1}g$ belongs to $U \Sigma U$. The fact that $\Sigma$ generates $G$ ensures that $\Gamma$ is connected. Moreover $G$ acts vertex-transitively by automorphisms on $\Gamma$. Since $U$ is compact open, the set  $U \Sigma U$ is a finite union of double cosets modulo $U$; this implies that $\Gamma$ is \textbf{locally finite}, i.e. the degree of each vertex is finite. Notice that all vertices have the same degree since $\Gamma$ is homogeneous. The graph $\Gamma$ is called a \textbf{Cayley--Abels graph} for $G$, since its construction was first envisaged by H.~Abels \cite[Beispiel~5.2]{Abels} and specializes to a Cayley graph when $G$ is discrete and $U = \{1\}$. Theorem~\ref{th:essentially_chief} is obtained by induction on the minimum degree of a Cayley--Abels graph.

 \subsection{Elementary groups}\label{sec:Elementary}

By its very nature, Theorem~\ref{th:essentially_chief} highlights the special role played by compact and discrete groups. A conceptual approach to studying the role of compact and discrete groups in the structure theory of tdlc groups is provided by P.~Wesolek's notion of \textbf{elementary groups}. That concept is inspired by the class of elementary amenable discrete groups introduced by M.~Day \cite{Day}. It is defined as the smallest class $\mathscr E$ of second countable tdlc groups (abbreviated \textbf{tdlcsc}) containing all countable discrete groups and all compact tdlcsc groups, which is stable under the following two group theoretic operations: 
\begin{itemize}
\item Given a tdlcsc group $G$ and a closed normal subgroup $N$, if $N \in \mathscr E$ and $G/N \in \mathscr E$, then $G \in \mathscr E$. In other words  $\mathscr E$ is \textit{stable under group extensions}. 
\item Given a tdlcsc group $G$ and a directed set $(O_i)_{i \in I}$ of open subgroups,  if $O_i\in \mathscr E$ for all $i$ and if   $G = \bigcup_i O_i$, then $G \in \mathscr E$. In other words  $\mathscr E$ is \textit{stable under directed unions of open subgroups}. 
\end{itemize}
(The class $\mathscr E$ has a natural extension beyond the second countable case, see \cite[\S6]{CRW_dense}. For simplicity of the exposition, we stick to the second countable case here.) 
Using the permanence properties of the class $\mathscr E$, it can be shown that every tdlcsc group $G$ has a largest normal closed subgroup that is elementary; it is denoted by $R_{\mathscr E}(G)$ and called the \textbf{elementary radical} of $G$. It indeed behaves as a radical, in the sense that it contains all elementary closed normal subgroups, and satisfies  $R_{\mathscr E}(G/R_{\mathscr E}(G)) = \{1\}$, see \cite[\S7.2]{Wes_elementary}. Further properties of the quotient $G/R_{\mathscr E}(G)$ will be mentioned in Section~\ref{sec:simple} below.

Similarly as for elementary amenable discrete groups, the class $\mathscr E$ admits a canonical rank function $\xi \colon \mathscr E \to \omega_1$, taking values in the set $\omega_1$ of countable ordinals,  called the \textbf{decomposition rank}. It measures the complexity of a given group $G \in \mathscr E$. By convention, the function $\xi$ is extended to all tdlcsc groups by setting $\xi(G) = \omega_1$ for each non-elementary tdlcsc group $G$. We refer to \cite{Wes_elementary}, \cite{Wes_survey} and \cite[\S 5]{Reid_survey}. Let us merely mention here that the class $\mathscr E$ has remarkable permanence properties (e.g. it is stable under passing to closed subgroups and quotient groups), that the rank function has natural monotonicity properties, and that a non-trivial compactly generated  group $G \in \mathscr E$ has a non-trivial discrete quotient. It follows in particular that if $G$ is a tdlcsc group having a closed subgroup $H \leq G$ admitting a non-discrete compactly generated topologically simple quotient, then  $G \not \in \mathscr E$. Therefore, the only compactly generated topologically simple groups in $\mathscr E$ are discrete. On the other hand, the class $\mathscr E$ contains numerous topologically simple groups that are not compactly generated, e.g. simple groups that are \textbf{regionally elliptic}, i.e. groups that can be written as a directed union of compact open subgroups. Those groups have decomposition rank~$2$.  Explicit examples appear in \cite[\S 3]{Willis_JAlg} or \cite[\S 6]{CapDM}.

 \subsection{More on chief factors}\label{sec:more-on-chief}
The existence of essentially chief series prompts us to ask whether the chief factors of $G$ depend upon the choice of a specific normal series in Theorem~\ref{th:essentially_chief}. 
It is tempting to tackle that question  by invoking arguments \`a la Jordan--H\"older. A technical obstruction for doing so is that the product of two closed normal subgroups need not be closed. More generally, given closed subgroups $A, N$ in $G$ such that $N$ is normal, the product $AN$ need not be closed so that the natural abstract isomorphism $ A/A\cap N\to AN/N$ need not be a homeomorphism. It is a continuous injective homomorphism of the locally compact group $A/A\cap N$ to a dense subgroup of the locally compact group $\overline{AN}/N$. This  illustrates the necessity of considering dense embeddings of locally compact groups. We shall come back to this theme in Section~\ref{sec:dense-and-local} below. In the context of chief factors, this has led Reid--Wesolek to define an equivalence relation on non-abelian chief factors of $G$, called \textbf{association}, defined as follows: the chief factors $K_1/L_1$ and $K_2/L_2$ are \textbf{associated} if $\overline{K_1L_2} = \overline{K_2L_1}$ and $K_i \cap \overline{L_1 L_2} = L_i$ for $i=1, 2$. In that case $K_1/L_1$ and $K_2/L_2$ both embed continuously as dense normal subgroups in $\overline{K_1L_2}/\overline{L_1L_2}$. We also recall that the \textbf{quasi-center} of a locally compact group $G$, denoted by $\QZ(G)$, is the collection of elements whose centralizer is open. It is a topologically characteristic (not necessarily closed) subgroup of $G$ containing all the discrete normal subgroups. It was first introduced by M.~Burger and S.~Mozes \cite{BuMo1}. Every non-trivial tdlcsc group with a dense quasi-center is elementary of decomposition rank~$2$ (see \cite[Lem.~5]{Reid_survey}).

\begin{thm}[{See \cite[Cor.~5]{Reid_survey}}]\label{th:invariance_chief}
Let $G$ be a compactly generated tdlc group and  let $\{1\} =A_0 < A_1 < A_2 < \dots < A_m = G$ and $\{1\} =B_0 < B_1 < B_2 < \dots < B_n = G$ be essentially chief series for $G$. Then for each $i \in \{0, 1, \dots, m\}$, if $A_i/A_{i-1}$ is a chief factor with a trivial quasi-center, there is a unique $j$ such that $B_j/B_{j-1}$ is a chief factor with a trivial quasi-center that is associated with $A_i/A_{i-1}$. In other words, the association relation establishes a bijection between the sets of chief factors with a trivial quasi-center appearing respectively in the two series. 
\end{thm}

The natural next question is to ask what can be said about chief factors. By the discussion above, one should focus on properties that are invariant under the association relation. Following Reid--Wesolek, an association class of non-abelian chief factors is called a \textbf{chief block}, and a group property shared by all members of a chief block is called a \textbf{block property}. The following are shown in  \cite{ReidWes_PLMS} to be block properties: compact generation, amenability, having a trivial quasi-center, having a dense quasi-center, being elementary of a given decomposition rank.

As mentioned above, every chief factor is topologically characteristically simple. In particular, a compactly generated chief factor is subjected to the following description. 

\begin{thm}[{See \cite[Cor.~D]{CaMo} and \cite[Rem.~3.10]{CaLB2}}]\label{thm:char-simple}
Let $G$ be a compactly generated non-discrete, non-compact tdlc group which is topologically characteristic simple. Then there is a compactly generated non-discrete topologically simple tdlc group $S$, an integer $d \geq 1$ and an injective continuous homomorphism   $S^d = S \times \dots \times S \to G$  of the direct product of $d$ copies of $S$, such that the image of each simple factor is a closed normal subgroup of $G$, and the image  of the whole product is dense. 
\end{thm}

In the setting of Theorem~\ref{thm:char-simple}, we say that $G$ is the \textbf{quasi-product} $d$ copies of the simple group $S$. 
Theorem~\ref{thm:char-simple} provides a major incentive to study the compactly generated non-discrete topologically simple tdlc groups. We shall come back to this theme in Section~\ref{sec:simple} below.  

Developing a meaningful structure theory for topologically characteristically simple tdlc groups that are not compactly generated is very challenging. Remarkably, significant results have been established by Reid--Wesolek \cite{ReidWes_PLMS} under the mild assumption of second countability (abbreviated sc). In spite  of the non-compact generation, they introduce an appropriate notion of chief blocks, and show that there are only three possible configurations for the arrangement of chief blocks in a   topologically characteristically simple tdlcsc group $G$, that they call \textbf{weak type}, \textbf{semisimple type}, and \textbf{stacking type}. Moreover, if $G$ is of weak type, then it is automatically elementary of decomposition rank~$\leq \omega +1$. The topologically characteristically simple groups $M(p)$ and $H(p)$ appearing in Examples~\ref{ex:1} and~\ref{ex:2} above are respectively of semisimple type and stacking type. We refer to \cite{ReidWes_PLMS} and \cite{Reid_survey} for details.

\section{Simple groups}\label{sec:simple}

Let $\mathscr S$ be the class of non-discrete, compactly generated, topologically simple locally compact groups and   $\mathscr S_{\mathrm{td}}$ be the subclass consisting of the totally disconnected members of $\mathscr S$. By the Gleason--Yamabe theorem \cite[Th.~6.0.11]{Tao}, all elements of $\mathscr S \setminus \mathscr S_{\mathrm{td}}$  are connected simple Lie groups. Prominent examples of groups in $\mathscr S_{\mathrm{td}}$ are provided by simple algebraic groups over non-Archimedean local fields, irreducible complete Kac--Moody groups over finite fields, certain groups acting on trees and many more, see \cite[Appendix~A]{CRW_part2}. A systematic study of the class $\mathscr S_{\mathrm{td}}$ as a whole has been initiated by Caprace--Reid--Willis in \cite{CRW_part2}, and continued with P.~Wesolek in \cite{CRW_dense} and with A.~Le Boudec in \cite{CaLB2}. We now outline some of their contributions. Another survey of the properties of non-discrete simple locally compact groups can be consulted in \cite{CapECM}; the present account emphasizes more recent results.
 
\subsection{Dense embeddings and local structure}\label{sec:dense-and-local}

As mentioned in Section~\ref{sec:more-on-chief} above, the failure of the second isomorphism theorem for topological groups naturally leads one to consider \textbf{dense embeddings}, i.e. continuous injective homomorphisms with dense image. 
If $G, H$ are locally compact groups and $\psi \colon H \to G$ is a dense embedding, and if $G$ is a connected simple Lie group or a simple algebraic group over a local field, then $H$ is discrete or $\psi$ is an isomorphism (see \cite[\S 3]{CRW_dense}). This property however generally fails for groups  $G \in \mathscr S$; see \cite{LB_CMH} for explicit examples. Nevertheless, as soon as the group $H$ is non-discrete, it turns out that key structural features of $G$ are inherited by the dense subgroup $H$. To state this more precisely, we recall  the  definition of the class $\mathscr R$ of robustly monolithic groups, introduced in \cite{CRW_dense}. A tdlc group $G$ is \textbf{robustly monolithic} if the intersection $M$ of all non-trivial closed normal subgroups of $G$ is non-trivial, if $M$ is topologically simple and if $M$ has a compactly generated  open subgroup without any non-trivial compact normal subgroup. The class $\mathscr R$ contains $\mathscr S_{\mathrm{td}}$ and that inclusion is strict. The following result provides the main motivation to enlarge one's viewpoint by considering $\mathscr R$ instead of the smaller class $\mathscr S_{\mathrm{td}}$. 

\begin{thm}[{See \cite[Th.~1.1.2]{CRW_dense}}]\label{thm:dense}
Let $G, H$ be tdlc groups and $\psi \colon H \to G$ be a dense embedding. If $G \in \mathscr R$ and $H$ is non-discrete, then $H \in \mathscr R$. 
\end{thm}

We emphasize that in general $H$ is not topologically simple even in the special case where $G \in \mathscr S_{\mathrm{td}}$. 

The approach in studying the classes $\mathscr S_{\mathrm{td}}$ and $\mathscr R$ initiated in  \cite{CRW_part2} is based on the concept of \textbf{locally normal subgroup}, defined as a subgroup whose normalizer is open. To motivate it, recall once more that if $M, N$ are closed normal subgroups of a tdlc group $G$, then the normal subgroup $MN$ need not be closed. On the other hand, if $U \leq G$ is a compact open subgroup, then $M \cap U$ and $N \cap U$ are closed normal subgroups of the compact group $U$ (hence they are both locally normal), so that the product $(M \cap U)(N \cap U)$ is closed. This observation motivates the definition of the \textbf{structure lattice} $\LN(G)$ of a tdlc group $G$, first introduced in \cite{CRW_part1}, defined as the set of closed locally normal subgroups of $G$, divided by the \textbf{local equivalence relation} $\sim$, where  $H \sim K$ if $H  \cap K$ is relatively open both in $H$ and in $K$. The local class of a closed locally normal subgroup $K$ is denoted by $[K]$. We also set $0 = [\{1\}]$ and $\infty = [G]$. The structure lattice carries a natural $G$-invariant order relation defined by the inclusion of representatives. The poset $\LN(G)$ is a modular lattice (see \cite[Lem.~2.3]{CRW_part1}). The greatest lower bound and least upper bound of two elements $\alpha, \beta \in \LN(G)$ are respectively denoted by $\alpha \wedge \beta$ and $\alpha \vee \beta$. When $G$ is a $p$-adic Lie group, the structure lattice $\LN(G)$ can naturally be identified with the lattice of ideals in the $\mathbf Q_p$-Lie algebra of $G$. The theory developed in \cite{CRW_part1} reveals that the structure lattice is especially well-behaved when the tdlc group $G$ is \textbf{[A]-semisimple}, i.e. $\QZ(G) = \{1\}$ and the only abelian locally normal subgroup of $G$ is $\{1\}$. That term is motivated by the fact that if $G$ is a $p$-adic Lie group, then it is [A]-semisimple if and only if $\QZ(G)= \{1\}$ and the $\mathbf Q_p$-Lie algebra of $G$ is semisimple, see \cite[Prop.~6.18]{CRW_part1}. An important result of P.~Wesolek is that the quotient $G/R_{\mathscr E}(G)$ of every tdlcsc group $G$ by its elementary radical is [A]-semisimple (see \cite[Cor.~9.15]{Wes_elementary}), so that every non-elementary group has a non-trivial [A]-semisimple quotient. The following result shows that [A]-semisimplicity is automatically fulfilled by groups in $\mathscr R$.  

\begin{thm}[{See \cite[Th.~A]{CRW_part2} and \cite[Th.~1.2.5]{CRW_dense}}]
Every group $G \in \mathscr R$ is [A]-semisimple. 
\end{thm}

Given an [A]-semisimple tdlc group $G$, two closed locally normal subgroups $H, K \leq G$ that are locally equivalent have the same centralizer; moreover they commute if and only if their intersection is trivial (see \cite[Th.~3.19]{CRW_part1}). This ensures that the map $\LN(G) \to \LN(G) : [K] \mapsto [K]^\perp = [C_G(K)]$ is well-defined, and that $\alpha \wedge \alpha^\perp = 0$ for all $\alpha \in \LN(G)$. This allows one to define the \textbf{centralizer lattice} of $G$ by setting $\LC(G) = \{\alpha^\perp \mid \alpha \in \LN(G)\}$. If $G$ is [A]-semisimple, the centralizer lattice $\LC(G)$ is a Boolean algebra (see  \cite[Th.~II]{CRW_part1}). We denote its Stone dual by $\Omega_G$. Thus $\Omega_G$ is a totally disconnected compact space endowed with a canonical continuous $G$-action by homeomorphisms. In general, the $G$-action on $\Omega_G$ need not be faithful. Actually, if $\LC(G) = \{0, \infty\}$ then $\Omega_G$ is a singleton. This happens if and only if any two non-trivial closed locally normal subgroups of $G$ have a non-trivial intersection. The following result shows that the dynamics of the $G$-action on $\Omega_G$ has remarkable  features. 

\begin{thm}[{See \cite[Th.~J]{CRW_part2} and \cite[Th.~1.2.6]{CRW_dense}}]\label{thm:G-boundary}
Let $G \in \mathscr R$. Then the $G$-action on $\Omega_G$ is minimal, strongly proximal, and has a compressible open set. Moreover the $G$-action on $\Omega_G$ is faithful if and only if $\LC(G) \neq \{0, \infty\}$. 
\end{thm}
	
Recall that a compact $G$-space $X$ is called \textbf{minimal} if every $G$-orbit is dense. It is called \textbf{strongly proximal} if the closure of each $G$-orbit in the space of probability measures on $X$ contains a Dirac mass. A non-empty subset $\alpha$ of $X$ is called \textbf{compressible} if for every non-empty open subset $\beta \subseteq X$ there exists $g \in G$ with $g\alpha \subseteq \beta$. Obviously, if $X$ is a minimal strongly proximal compact $G$-space and if $G$ fixes a probability measure on $X$, then $X$ is a singleton. Therefore, the following consequence of Theorem~\ref{thm:G-boundary} is immediate. 

\begin{cor}\label{cor:amenable}
Let $G \in \mathscr R$. If $G$ is amenable, then $\LC(G) = \{0,\infty\}$.
\end{cor}

A \textbf{local isomorphism} between tdlc groups $G_1, G_2$ is a triple $(\varphi, U_1, U_2)$ where $U_i$ is an open subgroup of $G_i$ and $\varphi \colon U_1 \to U_2$ is an isomorphism of topological groups. We emphasize that the structure lattice and the centralizer lattice are \textbf{local invariants}: they only depend on the local isomorphism class of the ambient tdlc group. However, for a group $G \in \mathscr R$, the compact $G$-space $\Omega_G$ can also be characterized by global properties among all compact $G$-spaces. In order to be more precise, let us first recall some terminology. Given an action of a group $G$ by homeomorphisms on a Hausdorff topological space $X$, we define the \textbf{rigid stabilizer} $\rist_G(U)$ of a subset $U \subseteq X$ as the pointwise stabilizer of the complement of $U$ in $X$. The $G$-action on $X$  is called \textbf{micro-supported} if for every non-empty open subset $U \subset X$ with $U \neq X$, the rigid stabilizer $\rist_G(U)$ acts non-trivially on $X$. The term ‘micro-supported’ was first coined in \cite{CRW_part2}, although the notion it designates has frequently appeared in earlier references, notably in the work of M.~Rubin on reconstruction theorems (see \cite{Rubin} and references therein). A prototypical example of a micro-supported of a tdlc group  is given by the action of the full automorphism group $\Aut(T)$ of a locally finite regular tree $T$ of degree~$\geq 3$ on the compact space $\partial T$ consisting of the ends of $T$. The following result shows that for a general group $G \in \mathscr R$, the  $G$-action on $\Omega_G$ shares many dynamical properties with the $\Aut(T)$-action on $\partial T$. 

\begin{thm}[{See \cite[Th.~J]{CRW_part2},   \cite[Th.~7.3.3]{CRW_dense} and \cite[Th.~7.5]{CaLB2}}]\label{thm:universal}
Let $G \in \mathscr R$. Then the $G$-action on $\Omega_G$ is micro-supported. Moreover, for each non-empty micro-supported compact $G$-space $X$ on which the $G$-action is faithful, there is a $G$-equivariant continuous surjective map $\Omega_G \to X$. In particular the $G$-action on $X$ is minimal, strongly proximal, and has a compressible open set. 
\end{thm}

This shows that $\Omega_G$ is universal among the faithful micro-supported compact $G$-spaces; in particular, the purely local condition that $\LC(G) = \{0, \infty\}$ ensures that $G$ does not have any faithful micro-supported continuous action on any compact space. Theorem~\ref{thm:universal} was first established for totally disconnected compact $G$-spaces in \cite{CRW_part2, CRW_dense}, and then extended to all compact $G$-spaces in \cite {CaLB2}, using tools from topological dynamics. Further properties of the $G$-space $\Omega_G$ and on the algebraic structure of groups in $\mathscr R$ can be consulted in those references. 

We now present another aspect of the local approach to the structure of simple tdlc groups. We define the \textbf{local prime content} of a tdlc group $G$, denoted by $\pi(G)$, to be the set of those primes $p$ such that every compact open subgroup $U \leq G$ contains an infinite pro-$p$ subgroup. 

\begin{thm}[{See \cite[Th.~H]{CRW_part2} and \cite[Cor.~1.1.4 and Th.~1.2.1]{CRW_dense}}]\label{th:prime}
The following assertions hold for any group $G \in \mathscr R$. 
\begin{enumerate}[label=(\roman*)]
\item The local prime content $\pi(G)$ is finite and nonempty. 
\item For each $p \in \pi(G)$, there is a group $G_{(p)} \in \mathscr R$ that is locally isomorphic to a pro-$p$ group, and a dense embedding $G_{(p)} \to G$. 
\item If $H$ is a tdlc group acting continuously and faithfully by automorphisms on $G$, then $H$ is locally isomorphic to a pro-$\pi(G)$ group. 
\end{enumerate}
\end{thm}
	
Roughly speaking, Theorem~\ref{th:prime}(ii) asserts that every group in $\mathscr R$ can be `approximated' by a locally pro-$p$ group in $\mathscr R$. The restriction on the automorphism group of a group in $\mathscr R$ from Theorem~\ref{th:prime}(iii) should be compared with the automorphism group of the  restricted product $M(p)$ from Example~\ref{ex:1}. Indeed, the Polish group $\mathrm{Sym}(\mathbf Z)$ embeds continuously in $\Aut(M(p))$ by permuting the simple factors, and every tdlcsc group continuously embeds in $\mathrm{Sym}(\mathbf Z)$. In some sense, the construction of stacking type chief factors in Example~\ref{ex:2} crucially relies on the hugeness of the group $\Aut(M(p))$. Theorem~\ref{th:prime}(iii) shows that the automorphism group of a group in $\mathscr R$ is considerably smaller. 

\bigskip
Let us finish this subsection with a brief discussion of classification problems. The work of S.~Smith \cite{Smith_Duke} shows that $\mathscr S_{\mathrm{td}}$ contains uncountably many isomorphism classes; his methods of proof  suggest that the isomorphism relation on $\mathscr S_{\mathrm{td}}$ has a similar complexity as the isomorphism relation on the class of finitely generated discrete simple groups. This provides evidence that the problem of  classifying groups in $\mathscr S_{\mathrm{td}}$ up to isomorphism is ill-posed. The recent results on the local structure of groups in $\mathscr S_{\mathrm{td}}$ or in $\mathscr R$ may be viewed as a hint to the fact the local isomorphism relation might be better behaved (see \cite[Th.~1.1.5]{CRW_dense}). At the time of this writing, we  do not know whether or not the groups in $\mathscr S_{\mathrm{td}}$ fall into countably many local isomorphism classes. However, classifying simple groups up to isomorphism remains a pertinent problem for some significant subclasses of $ \mathscr S_{\mathrm{td}}$. To wit, let us mention that, by \cite[Cor.~1.4]{CapStu}, a group $G \in \mathscr S_{\mathrm{td}}$ is isomorphic to a simple algebraic group over a local field if and only if it is locally isomorphic to a \textbf{linear group}, i.e. a subgroup of $\mathrm{GL}_d(k)$ for some integer $d$ and some locally compact field $k$. Lastly, a remarkable classification theorem concerning an important class of non-linear simple groups acting on locally finite trees  has been obtained by N.~Radu \cite{Radu}. 
It would be highly interesting to extend Radu's results by classifying all groups in $ \mathscr S_{\mathrm{td}}$ acting properly and continuously by automorphisms on a given locally finite tree $T$ in such a way that the action on the set of ends of $T$ is doubly transitive. That class is denoted by $ \mathscr S_T$. Results from \cite{CapRadu} ensure that the isomorphism relation restricted to $ \mathscr S_T$  is \textit{smooth} (see \cite[Definition 5.4.1]{Gao}), which means that it comes at the bottom of the hierarchy of complexity of classification problems in the formalism established by invariant descriptive set theory (see \cite[Ch.~15]{Gao}). Let us close this discussion by mentioning that we do not know whether   there is a tree $T$ such that $ \mathscr S_T$ contains uncountably many isomorphism classes.

\subsection{Applications to lattices}

 The study of lattices in semisimple Lie and algebraic groups has known tremendous developments since the mid 20th Century, with Margulis' seminal contributions as cornerstones. Remarkably, several key results on lattices have been established at a high level of generality, well  beyond the realm of linear groups. An early illustration  is provided by \cite{BernsteinKazhdan}. More recently, Y.~Shalom \cite{Shalom} and Bader--Shalom \cite{BaSha} have established an extension of Margulis' Normal Subgroup Theorem valid for all irreducible cocompact lattices in products of groups in $\mathscr S$, while various analogues of Margulis' superrigidity for irreducible lattices in products have been established for various kinds of target spaces, see \cite{Sha00, MonSha, Monod_JAMS, GKM, BaFu12, BaFu, CFI, Fioravanti}. Those results have in common that they rely on \textbf{transcendental methods}: they use a mix of tools from ergodic theory, probability theory and abstract harmonic analysis, but do not require any detailed consideration of the algebraic structure of the ambient group. Another breakthrough in this field was accomplished by M.~Burger and S.~Mozes \cite{BuMo2}, who constructed a broad family of new finitely presented infinite simple groups as irreducible lattices in products of non-linear groups in $\mathscr S_{\mathrm{td}}$. Their seminal work involves a mix of transcendental methods together with a fair amount of structure theory developed in \cite{BuMo1}. 

The following two results  rely in an essential way on the properties of the class $\mathscr S_{\mathrm{td}}$ outlined above. 

\begin{thm}[{See  \cite[Th.~A]{CaLB1}}]\label{thm:Wang1}
	Let $n \geq 2$ be an integer, let $G_1,\ldots,G_n \in \mathscr S_{\mathrm{td}}$ and  $ \Gamma \leq G = G_1 \times \dots \times G_n$ be a   lattice such that the projection $p_i(\Gamma)$ is dense in $G_i$ for all $i$. Assume that $\Gamma$ is cocompact, or that  $G$ has Kazhdan's property (T).  	Then the set of discrete subgroups of $G$ containing $\Gamma$ is finite. 
\end{thm}

\begin{thm}[{See  \cite[Th.~C]{CaLB1}}]\label{thm:Wang2} 
Let $n \geq 2$ be an integer and let $G_1,\ldots,G_n \in \mathscr S_{\mathrm{td}}$ be compactly presented. For every compact subset $K \subset G = G_1 \times \dots \times G_n$, the set of  discrete subgroups $\Gamma \leq G$ with $G = K\Gamma$ and with $p_i(\Gamma)$ dense in $G_i$ for all $i$,  is contained in a union of finitely many $\Aut(G)$-orbits.
\end{thm}

For a detailed discussion of the notion of \textbf{compactly presented} locally compact groups, we refer to \cite[Ch.~8]{CorHar}. 

Theorems~\ref{thm:Wang1} and~\ref{thm:Wang2} can be viewed as respective analogues of two theorems of H.~C. Wang \cite{Wang1, Wang2} on lattices in semisimple Lie groups and reveal the existence of positive lower bounds on the covolume of certain families of irreducible cocompact lattices. It should be underlined that the corresponding statements fail for lattices in a single group $G \in \mathscr S_{\mathrm{td}}$, see \cite[Th.~7.1]{BK90}. Theorem~\ref{thm:Wang2} is established by combining Theorem~\ref{thm:Wang1} with recent results on local rigidity of cocompact lattices in arbitrary groups, due to Gelander--Levit \cite{GeLe}. 

\subsection{Applications to commensurated subgroups}

Structure theory of tdlc groups provides valuable tools in exploring the so-called commensurated subgroups of an abstract group. In this section, we recall that connection and illustrate it with several recent results. Further results on commensurated subgroups will be mentioned in Section~\ref{sec:scale} below.

Let $\Gamma$ be a group. Two subgroups $\Lambda_1, \Lambda_2 \leq \Gamma$ are called \textbf{commensurate} if their intersection $\Lambda_1 \cap \Lambda_2$ has finite index both in $\Lambda_1$ and in $\Lambda_2$. The \textbf{commensurator} of a subgroup $\Lambda \leq \Gamma$, denoted by $\Comm_\Gamma(\Lambda)$, is the set of those $\gamma \in \Gamma$ such that $\Lambda$ and $\gamma \Lambda \gamma^{-1}$ are commensurate. It is easy to see that $\Comm_\Gamma(\Lambda)$ is a subgroup of $\Gamma$ containing the normalizer $N_\Gamma(\Lambda)$. The commensurator  has naturally appeared in group theory; one of its early occurrences is in Mackey's irreducibility criterion for induced unitary representations (see \cite{Mackey}). It also appears in a celebrated characterization of arithmetic  lattices in semisimple groups due to Margulis \cite[Ch.~IX, Th.~(B)]{Margulis}. A \textbf{commensurated subgroup} of $\Gamma$ is a subgroup $\Lambda \leq \Gamma$ such that $\Comm_\Gamma(\Lambda) = \Gamma$. Clearly, every normal subgroup of $\Gamma$ is commensurated; more generally, every subgroup that is commensurate to a normal subgroup is commensurated. Those commensurated subgroups are considered as trivial. For example, finite subgroups and subgroups of finite index are always commensurated subgroups. It is however important to underline that commensurated subgroups are not all of this trivial form.
 Indeed, an easy but crucial observation is that compact open subgroups are always commensurated. In particular, in the simple group $\PSL_2(\mathbf Q_p)$, the subgroup $\PSL_2(\mathbf Z_p)$ (which is obviously not commensurate to any normal subgroup of $\PSL_2(\mathbf Q_p)$) is commensurated. 
 
Let us next remark  that if $U$ is a commensurated subgroup of a group $G$ and  $\varphi \colon \Gamma \to G$ is a group homomorphism, then $\varphi^{-1}(U)$ is a commensurated subgroup of $\Gamma$. This is the case in particular if $G$ is a tdlc group and $U \leq G$ is a compact open subgroup. A fundamental observation is that all commensurated subgroups of $\Gamma$ arise in this way. More, precisely, a subgroup $\Lambda \leq \Gamma$ is commensurated if and only if there is a tdlc group $G$, a compact open subgroup $U \leq G$, and a homomorphism $\varphi \colon \Lambda \to G$ with dense image such that $\varphi^{-1}(U) = \Lambda$. Indeed, given a commensurated subgroup $\Lambda \leq \Gamma$, then $\Lambda$ acts on the coset space $\Gamma/\Lambda$ with finite orbits, so that the closure of the natural image of $\Gamma$ in the permutation group $\mathrm{Sym}(\Gamma/\Lambda)$, endowed with the topology of pointwise convergence, is a tdlc group containing the closure of the image of $\Lambda$ as a compact open subgroup. That tdlc group  is called the \textbf{Schlichting completion} of the pair $(\Gamma, \Lambda)$, denoted by $\Gamma/\!/\Lambda$. We refer to \cite{Schlichting},  \cite[Section 3]{ShaWil} and  \cite{ReidWesolek_Forum19} for more information. Let us merely mention here that a commensurated subgroup $\Lambda \leq \Gamma$ is commensurate to a normal subgroup if and only if the Schlichting completion $G = \Gamma/\!/\Lambda$ is \textbf{compact-by-discrete}, i.e. $G$ has a compact open normal subgroup (see \cite[Lem.~ 5.1]{CaLB2}). 

The occurrence of non-trivial commensurated subgroups in finitely generated  groups with few normal subgroups (e.g. simple groups, or \textbf{just-infinite groups}, i.e. groups all of whose proper quotients are finite) remains an intriguing phenomenon. On the empirical basis of the known examples, it seems to be rather rare. The following result provides valuable information  in that context.

\begin{thm}[{See \cite[Th.~5.4]{CaLB2}}]\label{thm:MicroComm}
Let $\Gamma$ be a finitely generated group. Assume that all normal subgroups of $\Gamma$ are finitely generated, and that every proper quotient of $\Gamma$ is virtually nilpotent. Let also $X$ be a compact $\Gamma$-space on which the $\Gamma$-action is faithful, minimal and micro-supported. Assume that at least one of the following conditions is satisfied:
\begin{enumerate}[label=(\arabic*)]
\item $\Gamma$ is residually finite.
\item $\Gamma$ fixes a probability measure on $X$.
\end{enumerate}

Then every commensurated subgroup of $\Gamma$ is commensurate to a normal subgroup.
\end{thm}

This applies to all finitely generated branch groups, as well as to numerous finitely generated almost simple groups arising in Cantor dynamics, and whose study has known spectacular recent developments (see \cite{Cor, Nekra} and references therein). We refer to \cite{CaLB2} for details and a more precise description of those applications. 

Let us briefly outline how the proof of Theorem~\ref{thm:MicroComm} works in the case where $\Gamma$ fixes a probability measure on $X$. Let $\Lambda \leq \Gamma$ be a commensurated subgroup and $G = \Gamma/\!/\Lambda$ be the corresponding Schlichting completion. That $\Gamma$ is finitely generated implies that $G$ is compactly generated. The hypotheses made on the normal subgroup structure of $\Gamma$ yield some restrictions on the essentially chief series of $G$ afforded by Theorem~\ref{th:essentially_chief}. More precisely, assuming by contradiction that $\Lambda$ is not commensurate to a normal subgroup, then the upper most chief factor $K/L$ with trivial quasi-center in an essentially chief series for $G$ must be compactly generated. Its structure is therefore described to by Theorem~\ref{thm:char-simple}. A key point in the proof, relying on various ingredients from topological dynamics and involving detailed considerations of the Chabauty space of closed subgroups of $\Gamma$ and $G$, is to show   that the given $\Gamma$-action on $X$ gives rise to a continuous, faithful, micro-supported $G/L$-action on a compact space $Y$ which is closely related to the original space $X$. Invoking (a suitable version of) Theorem~\ref{thm:universal} for the chief factor $K/L$ ensures that $Y$ has a compressible open set, from which it follows that $X$ has a compressible open set for the $\Gamma$-action. This finally contradicts the hypothesis of existence of a $\Gamma$-invariant probability measure.

\section{Scale methods}\label{sec:scale}

The scale of a tdlc group endomorphism, $\alpha$, is a positive integer that conveys information about the dynamics of the action of $\alpha$. Roughly speaking, $\alpha$ contracts towards the identity on one subgroup of $G$ and expands on another, and the scale is the expansion factor. This section gives an account of properties of the scale and descriptions of the action of $\alpha$ on certain associated subgroups of $G$ which, when applied to inner automorphisms, answer questions about group structure. 

Let $\alpha :G\to G$ be a continuous endomorphism. The \textbf{scale} of $\alpha$ is 
$$
s(\alpha) = \min\left\{ [\alpha(U) : \alpha(U)\cap U] \mid U\leq G\text{ compact and open}\right\}.
$$ 
This value is a positive integer, and the minimum is attained at subgroups that are said to be \textbf{minimizing for $\alpha$}, because $\alpha(U)\cap U$ is an open subgroup of the compact group $\alpha(U)$. The following results from \cite{Structure,Further,Endo} relate minimising subgroups to the dynamics of $\alpha$.

\begin{thm}\label{thm:tidy}
Let $\alpha$ be a continuous endomorphism of the tdlc group $G$ and let $U\leq G$ be compact and open. Define subgroups
\begin{align*}
U_+ &= \left\{ u\in U \mid \exists \{u_n\}_{n\geq0}\subset U\text{ with }u_0 = u\text{ and }u_n = \alpha(u_{n+1})\right\}\\
U_- &= \left\{ u\in U \mid \alpha^n(u)\in U\text{ for all }n\geq0\right\}.
\end{align*}
Also define the subgroup $U_{--} = \bigcup_{n\geq0} \alpha^{-n}(U_-)$ of $G$. \\
Then $U$ is minimizing for $\alpha$ if and only if \\
\phantom{a}\hfil {\textbf{TA}}: $U = U_+U_-$\quad  and \quad {\textbf{TB}}: $U_{--}$ is closed.\hfil
\end{thm}
A compact open subgroup $U$ satisfying {\textbf{TA}} and {\textbf{TB}} is said to be \textbf{tidy for $\alpha$}, and $s(\alpha) = [\alpha(U_+) : U_+]$ for any such subgroup $U$. Tidiness has two further dynamical interpretations: (1) an $\alpha$-trajectory $\{\alpha^n(g)\}_{n\geq0}$ cannot return to a tidy subgroup once it departs; and (2) when $\alpha$ is an automorphism, $U$ is tidy for $\alpha$ if and only if the orbit $\{\alpha^n(U)\}_{n\in\mathbb{Z}}$ is a geodesic for the metric $d(U,V) = \log[U:U\cap V] + \log[V:U\cap V]$ on the set of compact open subgroups of $G$. 

Note that every compact open subgroup of $G$ has a subgroup $U$ for which \textbf{TA} holds and, if $\alpha$ is the inner automorphism $\alpha_g(x) := gxg^{-1}$, then property \textbf{TA} implies that $Ug^mUg^nU = Ug^{m+n}U$ for all $m,n\geq0$. These points were already used in \cite{Bernshtein} in the proof that a reductive group over a locally compact field of positive characteristic is type I, where they were observed to hold in such groups.

In the following compilation of results from \cite{Structure,Moller,Further,Endo}, $\Delta$ denotes the modular function on the automorphism group of $G$.
\begin{thm}\label{thm:scale}
The scale $s : \End(G) \to \mathbb{Z}^+$ satisfies: 
\begin{enumerate}[label=(\roman*)]
\item $s(\alpha)=1$ if and only if there is a compact open subgroup $U\leq G$ with $\alpha(U)\leq U$;
\item  $s(\alpha) = \lim_{n\to\infty} [\alpha^n(V) : \alpha^n(V)\cap V]^{\frac{1}{n}}$ for every compact open $V\leq G$ and $s(\alpha^n)=s(\alpha)^n$ for every $n\geq0$; and 
\item if $\alpha$ is an automorphism, then $\Delta(\alpha) = s(\alpha)/s(\alpha^{-1})$.
\end{enumerate}
The function $s\circ \alpha_{{\bullet}} : G\to \mathbb{Z}^+$, with $\alpha_g(x) = gxg^{-1}$, is continuous for the group topology on $G$ and the discrete topology on $\mathbb{Z}^+$.
\end{thm}
Continuity of $s\circ \alpha_{{\bullet}}$ is implied by the fact that, if $U$ is tidy for $g$, then $U$ is also tidy for all $h\in UgU$ and $s(h) = s(g)$, \cite[Theorem 3]{Structure}. 

Questions about the structure of tdlc groups may be answered with scale and tidy subgroup techniques. K.~H.~Hofmann and A.~Mukherjea conjectured in \cite{HofMuk} that all locally compact groups are `neat' -- a property involving the conjugation action by a single element $g$. They used approximation by Lie groups to reduce to the totally disconnected case, and subgroups tidy for~$g$ are used in \cite{JaRoWi} to show that all groups are neat. Answering another question of K.~H.~Hofmann, the set ${\sf per}(G)$ comprising those elements of $G$ such that the closure of $\langle g\rangle$ is compact is shown in \cite{Periodic} to be closed by appealing to the properties of the scale given in Theorem~\ref{thm:scale}.

The scale and the subgroup $U_+$ associated with it in Theorem~\ref{thm:tidy} are given a concrete representation in \cite{ContractionB}. Put $U_{++} = \bigcup_{n\geq0} \alpha^n(U_+)$. Then $U_{++}$ is closed if $U$ is tidy and $U_{++}\rtimes \langle\alpha\rangle$ acts on a regular tree with valency $s(\alpha)+1$: the image of $U_{++}\rtimes \langle\alpha\rangle$ is a closed subgroup of the isometry group of the tree; is transitive on vertices; and fixes an end of the tree. The resulting isometry groups of trees correspond to the \textit{self-replicating groups} studied in \cite{Nekrashevych}. Moreover, the semi-direct product $U_{++}\rtimes \langle\alpha\rangle$ also  belongs to the family of \textit{focal hyperbolic groups} studied in \cite{CCMT}.

\subsection{Contraction and other groups}
\label{sec:contraction}

Subgroups of $G$ defined in terms of the action of $\alpha$ are related to the scale and tidy subgroups. It is convenient to confine the statements to automorphisms here. Extensions to endomorphisms may be found in \cite{Endo,BGTEndo}.

The \textbf{contraction subgroup} for $\alpha\in \Aut(G)$ is 
$$
\con(\alpha) = \left\{ x\in G \mid \alpha^n(x)\to\ident \text{ as }n\to\infty\right\}.
$$ 
The next result, from \cite{ContractionB, Jaworski}, relates contraction subgroups to the scale.
\begin{thm}
Let $\alpha\in\Aut(G)$. Then $\bigcap\left\{ U_{--}\mid U\text{ is tidy for }\alpha\right\}$ is equal to $\overline{\con(\alpha)}$, and $s(\alpha^{-1})$ is equal to the scale of the restriction of $\alpha^{-1}$ to $\overline{\con(\alpha)}$. Hence $s(\alpha^{-1})>1$ if and only if $\overline{\con(\alpha)}$ is not compact.
\end{thm}
If $G$ is a $p$-adic Lie group, then $\con(\alpha)$ is closed for every $\alpha$, \cite{WangJSP}, but that is not the case if, for example, $G$ is the isometry group of a regular tree, or a certain type of complete Kac-Moody group \cite{contractionKM}, or if $\LC(G) \ne \{0,\infty\}$ \cite{CRW_part2}. Closedness of $\con(\alpha)$ is equivalent, by \cite[Theorem 3.32]{ContractionB}, to triviality of the \textbf{nub} subgroup, 
$$
\nub(\alpha) = \bigcap\left\{ U\mid U\text{ tidy for }\alpha\right\}.
$$ 
The nub for $\alpha$ is compact and is the largest $\alpha$-stable subgroup of $G$ on which $\alpha$ acts ergodically, which sharpens the theorem of N.~Aoki in \cite{Aoki} that a totally disconnected locally compact group with an ergodic automorphism must be compact. P.~Halmos had asked in \cite{Halmos} whether that was so for all locally compact groups. See \cite{KaufmanRajagopalan}, \cite{WuHalmos} for the connected case, and also \cite{PrevitsWu}.

The structure of closed contraction subgroups $\con(\alpha)$ is described precisely in \cite{ContractionG}. If $\con(\alpha)$ is closed, there is a composition series
$$
\{\ident\} =  G_0 \triangleleft \cdots \triangleleft G_n = \con(\alpha)
$$
of $\alpha$-stable closed subgroups of $\con(\alpha)$ such that the factors $G_{i+1}/G_i$ have no proper, non-trivial $\alpha$-stable closed subgroups. The factors appearing in any such series are unique up to permutation and isomorphism, and their isomorphism types come from a countable list: each torsion factor being a restricted product $\bigoplus_{i \in \mathbf{Z}} (G_i, U_i)$ with $G_i=F$, a finite simple group, and $U_i=F$ if $i\geq0$ and trivial if $i<0$, and the automorphism the shift; and each divisible factor being a $p$-adic vector group and the automorphism a linear transformation. Moreover, $\con(\alpha)$ is the direct product $T\times D$ with $T$ a torsion and $D$ a divisible $\alpha$-stable subgroup. The divisible subgroup $D$ is a direct product $D_{p_1}\times \dots \times D_{p_r}$ with $D_{p_i}$ a nilpotent $p$-adic Lie group for each $p_i$. The torsion group $T$ may include non-abelian irreducible factors but, should it happen to be locally pro-$p$, then it is nilpotent too, see \cite{ContractionG3}. The number of non-isomorphic locally pro-$p$ closed contraction groups is uncountable, \cite{ContractionG2}.

Contraction groups correspond to unipotent subgroups of algebraic groups and, following \cite{TitsCore}, the \textbf{Tits core}, $G^\dag$, of the tdlc group $G$ is defined to be the subgroup generated by all closures of contraction groups. It is shown in \cite{CRW} that, if $G$ is topologically simple, then $G^\dag$ is either trivial or is abstractly simple and dense in $G$.  

The correspondence with algebraic groups is pursued in \cite{ContractionB}, where the \textbf{parabolic subgroup} for $\alpha\in\Aut(G)$ is defined to be
$$
\parb(\alpha) = \left\{ x\in G \mid \{\alpha^n(x)\}_{n\geq0} \text{ has compact closure}\right\},
$$  
and the \textbf{Levi factor} to be $\lev(\alpha) = \parb(\alpha)\cap \parb(\alpha^{-1})$. Then $\parb(\alpha)$, and hence $\lev(\alpha)$, is closed in $G$, \cite[Proposition 3]{Structure}. It may be verified that $\con(\alpha)\triangleleft \parb(\alpha)$ and shown, see \cite{ContractionB}, that $\parb(\alpha) = \lev(\alpha)\con(\alpha)$.

\subsection{Flat groups of automorphisms}
\label{sec:flat}

A group, $\mathcal{H}$, of automorphisms of $G$ is \textbf{flat} if there is a compact open subgroup, $U\leq G$, that is tidy for every $\alpha\in\mathcal{H}$. The stabilizer of $U$ in $\mathcal{H}$ is called the \textbf{uniscalar subgroup} and denoted $\mathcal{H}_u$. The factoring of subgroups tidy for a single automorphism in Theorem~\ref{thm:tidy} extends to flat groups as follows.
\begin{thm}[\cite{SimulTriang}]\label{thm:flat1}
Let $\mathcal{H}$ be a finitely generated flat group of automorphisms of $G$ and suppose that $U$ is tidy for $\mathcal{H}$. Then $\mathcal{H}_u\triangleleft \mathcal{H}$ and there is $r\geq0$ such that
$$
\mathcal{H}/\mathcal{H}_u \cong \mathbb{Z}^r.
$$
\begin{itemize}
\item There are $q\geq0$ and closed groups $U_j\leq U$, $j\in \{0,1,\dots, q\}$ such that: $\alpha(U_0) = U_0$; $\alpha(U_j)$ is either a subgroup or supergroup of $U_j$ for every $j\in\{1,\dots,q\}$; and $U=U_0U_1\dots U_q$.
\item $\widetilde{U}_j := \bigcup_{\alpha\in\mathcal{H}} \alpha(U_j)$ is a closed subgroup of $G$ for each $j\in\{1,\dots,q\}$. 
\item There are, for each $j\in\{1,\dots,q\}$, an integer $s_j>1$ and a surjective homomorphism $\rho_j : \mathcal{H}\to (\mathbb{Z},+)$ such that $\Delta(\alpha|_{\widetilde{U}_j}) = s_j^{\rho_j(\alpha)}$. 
\item The integers $r$ and $q$, and integers $s_j$ and homomorphisms $\rho_j$ for each $j\in\{1,\dots,q\}$, are independent of the subgroup $U$ tidy for $\mathcal{H}$.
\end{itemize}
\end{thm}
The number $r$ in Theorem~\ref{thm:flat1} is the \textbf{flat rank} of $\mathcal{H}$. The singly-generated group $\langle \alpha\rangle$ has flat rank equal to $0$ if $\alpha$ is uniscalar and $1$ if not. Flat groups of automorphisms with rank at least $1$ correspond to Cartan subgroups in Lie groups over local fields and to stabilisers of apartments in isometry groups of buildings, \cite{contractionBRW}. 

More generally, flatness of groups of automorphisms may be shown by the following converse to the fact that flat groups are abelian modulo the stabliser of tidy subgroups. 
\begin{thm}[\cite{SimulTriang}\cite{ShaWil}]\label{thm:flat2}
Every finitely generated nilpotent subgroup of $\Aut(G)$ is flat, and every polycyclic subgroup is virtually flat.
\end{thm}

Flatness is used -- in combination with: bounded generation of arithmetic groups, \cite{Tavgen, Witte-Morris}; the fact that almost normal subgroups are close to normal, \cite{BergmanLenstra}; and the Margulis normal subgroup theorem, \cite{Margulis} -- to prove the Margulis-Zimmer conjecture in the special case of Chevalley groups in \cite{ShaWil} and show that there are no commensurated subgroups of arithmetic subgroups other than the natural ones.

\section{Future directions}\label{sec:future}

The contributions to the structure  theory of tdlc groups surveyed in this article highlight that, for a general tdlc group $G$, as soon as the topology is non-discrete, its interaction with the group structure yields significant algebraic constraints. As mentioned in the introduction, we view the dynamics   of the conjugation action as a unifying theme of our considerations. 
The results we have surveyed reveal that those dynamics   tend to be richer than one might expect. This is especially the case among tdlc groups that are non-elementary.  We hope that further advances will shed more light on this paradigm in the future. 

Concerning decomposition theory, it is an important open problem to clarify what distinguishes elementary and non-elementary tdlc groups. A key question asks whether every non-elementary tdlcsc group $G$ contains a closed subgroup $H$ admitting a quotient in $\mathscr S_{\mathrm{td}}$. Concerning simple groups,  our results yield a dichotomy, depending on whether the centralizer lattice is trivial or not. The huge majority of known examples of groups in  $\mathscr S_{\mathrm{td}}$ (listed in \cite[Appendix~A]{CRW_part2}) have a non-trivial centralizer lattice, the most notable exceptions being the simple algebraic groups over local fields. Finding new groups in $\mathscr S_{\mathrm{td}}$ with a trivial centralizer lattice would be a decisive step forward. A fundamental source of examples of tdlc groups is provided by Galois groups of transcendental field extensions with finite transcendence degree (see \cite[Th.~2.9]{Rovinsky}, highlighting the occurrence of topologically simple groups), but this territory remains largely unexplored from the  viewpoint of structure theory of tdlc groups.  Concerning scale methods, the structure of tdlc groups all of whose elements are \textbf{uniscalar} (i.e. have scale~$1$) is still mysterious. In particular, we do not know whether every such group is elementary. This is equivalent to asking whether a tdlc group, all of whose closed subgroups are unimodular, is necessarily elementary. A positive answer would provide a formal incarnation to the claim that the dynamics of the conjugation action is non-trivial for all non-elementary tdlc groups. 
We refer to  \cite{CaMo_Future} for a more extensive list of specific problems. 

We believe that a good measurement of the maturity of a mathematical theory is provided by its ability to solve problems arising on the outside of the theory. For the structure theory of tdlc groups, the Margulis--Zimmer conjecture appears as a natural target. As mentioned in Section~\ref{sec:scale}, partial results in the non-uniform case, relying on scale methods on tdlc groups, have already been obtained in \cite{ShaWil}.

Another source of external problems is provided by abstract harmonic analysis. As mentioned in the introduction, the emergence of locally compact groups as an independent subject of study coincides with the foundation of abstract harmonic analysis. However, fundamental problems clarifying the links between the algebraic structure of a locally compact group and the properties of its  unitary representations remain open. The class of amenable locally compact groups is defined by a representation theoretic property (indeed, a locally compact group is amenable if and only if every unitary representation is weakly contained in the regular), but purely algebraic characterizations of amenable groups are still missing. In particular, the following non-discrete version of Day's problem is open and intriguing: \textit{Is every amenable  second countable tdlc group elementary (in the sense of Section~\ref{sec:dec})?}  The unitary representation theory also reveals a fundamental dichotomy between locally compact groups of \textbf{type I} (roughly speaking, those for which the problem of classifying the irreducible unitary representations up to equivalence is tractable) and the others (see \cite{BeHa, Dixmier, Mackey}). Algebraic characterizations of type I groups are also desirable. In particular, we underline the following question: \textit{Does every second countable locally compact group of type I contain a cocompact amenable subgroup?} For  a more detailed discussion of that problem and related results, we refer to \cite{CKM}. 

\subsection*{Acknowledgements}

It is a pleasure to take this opportunity to heartily thank our collaborators, past and present, for their ideas and companionship. We are grateful to  Adrien Le Boudec and Colin Reid for their comments on a preliminary version of this article.

\end{document}